\documentclass{easychair}
\usepackage{amstext}
\usepackage{amssymb}

\usepackage{doc}

\newtheorem{theorem}{Theorem}

\newtheorem{definition}{Definition}

\newtheorem{fact}{Fact}

\global\long\def\qftp{\operatorname{qftp}}%
\global\long\def\range{\operatorname{rng}}
\global\long\def\Aut{\operatorname{Aut}}%

\title{Intersecting sets in probability spaces and Shelah's classification}

\author{
Artem Chernikov\inst{1}\thanks{Chernikov was partially supported by the NSF
CAREER grant DMS-1651321 and by the NSF Research Grant DMS-2246598.}
\and
    Henry Towsner\inst{2}\thanks{Towsner was partially supported by the NSF grant DMS-2054379.}
}

\institute{
 University of Maryland, 
  College Park, MD, US and University of California, Lost Angeles, CA, US\\
  \email{artem@umd.edu}
  \and
   University of Pennsylvania,
 Philadelphia, PA, US\\
   \email{htowsner@math.upenn.edu}
 }

\authorrunning{Artem Chernikov}

\titlerunning{Intersecting sets in probability spaces and Shelah's classification}

\begin{document}

\maketitle

\begin{abstract}
For  $n \in \mathbb{N}$ and $\varepsilon > 0$, given a sufficiently long sequence of events in a probability space all of measure at least $\varepsilon$, some $n$ of them will have a common intersection. A more subtle pattern: for any $0 < p < q < 1$, we cannot find events $A_i$ and $B_i$ so that $\mu \left( A_i \cap B_j \right) \leq p$  and $\mu \left( A_j \cap B_i\right) \geq q$ for all $1 < i < j < n$, assuming $n$ is sufficiently large. This is closely connected to model-theoretic stability of probability algebras. We survey some results from our recent work in \cite{CheTow} on more complicated patterns that arise when our events are indexed by multiple indices. In particular, how such results are connected to higher arity generalizations of de Finetti's theorem in probability, structural Ramsey theory, hypergraph regularity in combinatorics, and model theory.
\end{abstract}



%
%

\section{Intersections in a sequence of sets of positive measure}\label{sec: classic fact}

The following is a basic fact on intersections of sets of positive measure in probability spaces (there exist more precise infinitary/density versions, e.g.~Bergelson's lemma in dynamics \cite{bergelson1985sets}):
\begin{fact}\label{fac: basic}
	For every $\varepsilon \in \mathbb{R}_{>0}$ and $n \in \mathbb{N}$ there exists  $N \in \mathbb{N}$ satisfying the following. If $(X, \mathcal{B}, \mu)$ is a probability space (i.e.~$\mathcal{B}$ is a $\sigma$-algebra of subsets of $X$ and $\mu$ is a countably additive probability measure on $\mathcal{B}$) and $A_i \in \mathcal{B}$ are measurable sets with $\mu(A_i) \geq \varepsilon$ for $1 \leq i \leq N$, then $\mu \left( \bigcap_{i \in I} A_i \right) > 0$ for some $I \subseteq [N] = \{1, \ldots, N \}$ with $|I|=n$.
\end{fact}
\noindent We sketch an overcomplicated proof of this fact  in the remainder of Section \ref{sec: classic fact}, as a warm up for what comes later. 
If the random variables $\textbf{1}_{A_i}: X \to \{0,1\}$ in Fact \ref{fac: basic}  were independent, then of course $\mu \left( \bigcap_{i \in [n] } A_i\right) = \prod_{i \in [n]} \mu(A_i) \geq \varepsilon^n > 0$.
We will reduce to this case. 
Assume from now on that for some fixed $\varepsilon >0$ and $n$, no $N \in \mathbb{N}$ satisfies the claim.

\subsection{Homogenizing the sequence}\label{sec: homogenize sequence}

 Using Ramsey's theorem, we can homogenize our sequence arbitrarily well. E.g., we could assume that for any fixed $\delta > 0$ and $k$,  $\mu(A_{i_1} \cap \ldots \cap A_{i_k}) \approx^{\delta} \mu(A_{j_1} \cap \ldots \cap A_{j_k})$ for any $i_1 < \ldots < i_k$, $j_1 < \ldots < j_k$, and similarly for  arbitrary Boolean combinations of the $A_i$'s.

 Using a compactness argument (e.g.~taking Loeb measure on an ultraproduct of counterexamples), we can thus find some large probability space $(X, \mathcal{B}, \mu)$ and sets $ A_i \in \mathcal{B}, \mu(A_i) \geq \varepsilon$ for $i \in \mathbb{N}$, still intersection of any $n$ of them has measure $0$, so that the sequence of random variables $\left(\textbf{1}_{A_i} : i \in \mathbb{N} \right)$ is \emph{spreadable}.

\subsection{de Finetti's theorem}\label{sec: de finetti}

\begin{definition}
	A sequence of $[0,1]$-valued random variables $(\xi_i : i \in \mathbb{N})$ is \emph{spreadable} if for every $n \in \mathbb{N}$ and $i_1 < \ldots < i_n, j_1 < \ldots < j_n$ we have $(\xi_{i_1}, \ldots, \xi_{i_n}) =^{\textrm{dist}} (\xi_{j_1}, \ldots, \xi_{j_n})$.
\end{definition}
For example, every i.i.d. (independent, identically distributed) sequence of random variables is spreadable. The converse holds ``up to mixing'':

\begin{fact}[de Finetti's theorem]
	If an infinite sequence of random variables $(\xi_i : i \in \mathbb{N})$ on $(X, \mathcal{B}, \mu)$ is spreadable then there exists a probability space  $(X',\mathcal{B}',\mu')$,  a Borel function $f: [0,1]^2 \to [0,1]$ and a collection of $\textrm{Uniform}[0,1]$ \emph{i.i.d.}~random variables $\left\{ \zeta_i : i \in \mathbb{N} \right \} \cup \{\zeta_{\emptyset}\}$ on $X'$ so that 
$\left( \xi_i : i \in \mathbb{N} \right)=^{\textrm{dist}} 
	\left( f \left(  \zeta_i, \zeta_{\emptyset} \right) : i \in \mathbb{N} \right)$.
\end{fact}
	 This gives us an i.i.d.~counterexample to Fact \ref{fac: basic}, and we can conclude.
\subsection{Exchangeable versus spreadable sequences}\label{sec: exch vs spread}
More precisely, de Finetti obtained this conclusion under a stronger assumption that the sequence $\left( \xi_i : i \in \mathbb{N} \right)$ is \emph{exchangeable}, that is:  for any $n \in \mathbb{N}$, any permutation $\sigma \in \textrm{Sym}(n)$ and $i_1 < \ldots < i_n$ we have $(\xi_{i_1}, \ldots, \xi_{i_n}) =^{\textrm{dist}} (\xi_{i_{\sigma(1)}}, \ldots, \xi_{j_{\sigma(n)}})$.
And then Ryll-Nardzewski \cite{ryll1957stationary} proved that exchangeability is equivalent to spreadability.
Curiously, Ryll-Nardzewski has a well-known theorem in model theory, but here he worked as a probabilist. It turns out that this result is closely connected to \emph{stability} --- a central notion in modern model theory.

\section{Model theoretic stability of probability algebras}

Modern model theory begins with \emph{Morley's Categoricity Theorem}: for a countable theory $T$, if it has only one model of some uncountable cardinality (up to isomorphism), then it has only one model of \emph{every} uncountable cardinality.
Morley  conjectured \cite{morley1965categoricity} a generalization: 	for a countable theory $T$, the number of its models of size $\kappa$ is non-decreasing on uncountable $\kappa$.

In his solution of Morley's conjecture \cite{shelah1990classification}, Shelah isolated the importance of  \emph{stable theories} and developed a lot of machinery to analyze  models of stable theories. Stability was rediscovered many times in various contexts, e.g.~by Grothendieck in his work on Banach spaces, in dynamics as WAP systems (Weakly Almost Periodic), in  machine learning as Littlestone dimension, etc.

In particular, probability algebras are stable, viewed as structures in continuous logic. This is implicit in Ryll-Nardzewski's theorem (``every indiscernible sequence is totally indiscernible''), later in Krivine and Maurey \cite{krivine1981espaces}, explicit in Ben Yaacov \cite{yaacov2013theories}. A more general version  was given by Hrushovski (proved using array de Finetti, discussed in Section \ref{sec: exch arrays}), and Tao gave a short elementary proof \cite{Tao}:
	\begin{fact}\label{fac: stability}
	For any $0 \leq p < q \leq 1$ there is $N$ satisfying: if  $(X,\mathcal{B},\mu)$ is a probability space, and $A_1, \ldots, A_n, B_1, \ldots, B_n \in \mathcal{B}$ satisfy $\mu(A_i \cap B_j) \geq q$ and $\mu(A_j \cap B_i) \leq p$ for all $1 \leq i < j \leq n$, then $n \leq N$.
\end{fact}
This result has many applications: Hrushovski's work on approximate subgroups \cite{hrushovski2012stable}, Tao's algebraic regularity lemma \cite{tao2015expanding}, work in topological dynamics by Tsankov, Ibarlucia \cite{ibarlucia2021model}, etc.

\section{Intersecting multi-parametric families of events}
We obtain a \emph{higher arity} generalization of Fact \ref{fac: basic}:
\begin{theorem}[Chernikov, Towsner  \cite{CheTow}]\label{thm: main1}
For every  finite bipartite graph $H = (V_0,W_0, E_0)$  and $\varepsilon \in (0,1]$ there exists a finite bipartite graph $G = (V,W, E)$ and $\delta > 0$ (depending only on $H$ and $\varepsilon$) satisfying the following. Assume that  $(X, \mathcal{B}, \mu)$ is a probability space, and for every $(v, w) \in V \times W$ a measurable set  $A_{v,w} \in \mathcal{B}$ so that: for any $(v,w) \in E, (v',w') \notin E$ we have $\mu(A_{v,w}) - \mu(A_{v',w'}) \geq \varepsilon$. Then there exists an induced subgraph $H' = (V', W', E')$  of $G$ (i.e.~$V' \subseteq V, W' \subseteq W$ and $E' = E \cap  (V' \times W')$) isomorphic to $H$ so that: 
$$\mu \left( \left(\bigcap_{(v,w) \in E'} A_{v,w} \right) \cap \left( \bigcap_{(v,w) \in (V' \times W') \setminus E'} X \setminus A_{v,w}  \right) \right) \geq \delta.$$
\end{theorem}

More precisely, Theorem \ref{thm: main1} follows from  \cite[Lemma 10.13]{CheTow} and compactness. With high probability, a sufficiently large $G$ taken at random will work. More generally, we prove it there for partite hypergraphs of any arity instead of just graphs.  The question is motivated by \emph{Keisler randomizations} of first-order structures \cite{keisler1999randomizing} and whether they preserve NIP (Ben Yaacov, related to work of Talagrand on VC-dimension for functions \cite{yaacov2009continuous}) and its higher arity generalization \emph{$n$-dependence} (where Ben Yaacov's analytic proof for $n=1$ does not seem to generalize).

In what follows we outline a proof of Theorem \ref{thm: main1}. The overall strategy is similar to the proof above for sequences of events, but each of the steps becomes harder.

\subsection{Structural Ramsey theory, and infinite limits of Ramsey classes}

Let $\mathcal{K}$ be a class of finite $\mathcal{L}_{0}$-structures, where $\mathcal{L}_{0}$ is a relational language (for example, finite graphs). For $A,B\in K$, let ${B \choose A}$ be the set of all $A'\subseteq B$ s.t.~$A'\cong A$ (we work with substructures instead of embeddings for simplicity).

\begin{definition}
	$\mathcal{K}$ is \emph{Ramsey} if for any $A,B\in K$ and $k\in\omega$ there is some $C\in K$ s.t. for any coloring $f:{C \choose A}\to k$, there is some $B'\in{C \choose B}$ s.t. $f\restriction{B' \choose A}$ is constant.
\end{definition}

 The usual Ramsey theorem means: the class of finite linear orders is Ramsey. The subject of structural Ramsey theory started with the following fundamental result of Nes\'etril, R\"odl  \cite{nevsetvril1977partitions} and Abramson, Harrington \cite{abramson1978models}:
\begin{fact}\label{fac: graphs Ramsey}
	For any $k\in \mathbb{N}_{\geq 1}$, the class of all finite ordered $k$-hypergraphs is Ramsey.
\end{fact}

\begin{fact}
	 Given a Ramsey class  $\mathcal{K}$ of finite structures, there exists a unique (up to isomorphism) countable structure $\widetilde{K}$ (called the \emph{Fra\"iss\'e limit} of $\mathcal{K}$) so that the class of its finite substructures is precisely $\mathcal{K}$ and $\widetilde{K}$ is \emph{homogeneous}, i.e.~if $K_0$ and $K_1$ are finite substructures of $\widetilde{K}$ and $f: K_0 \to K_1$ is an isomorphism, then $f$ extends to an automorphism of the whole structure $\widetilde{K}$.
\end{fact}
 E.g., if $\mathcal{K}$ is the class of all graphs, its limit $\widetilde{K}$ is the countable Rado's random graph; and if $\mathcal{K}$ is the class of finite linear orders, then its limit is $(\mathbb{Q}, <)$.

Understanding which structures are Ramsey is an active subject, with connections to model theory and topological dynamics (Ramsey property of $\mathcal{K}$ is equivalent to the extreme amenability of the  group $\Aut(\widetilde{K})$ --- via Kechris-Pestov-Todorcevic correspondence \cite{kechris2005fraisse}).

\subsection{Finding an ``exchangeable'' counterexample}\label{sec: exch counterex}

 For any $k\in \mathbb{N}_{\geq 1}$, using that the class of all finite ordered  partite $k$-hypergraphs is Ramsey (viewed as structures in the language $E, P_1, \ldots, P_k, <$ with $P_i$ a partition of vertices, $E \subseteq P_1 \times \ldots \times P_k$ and $P_i < P_j$ for $i<j$, e.g.~\cite[Appendix A]{chernikov2019n}), we let  $\mathcal{OH}_k$ denote its \emph{Fra\"iss\'e limit}. And we let $\mathcal{H}_k$ be its reduct forgetting the ordering.

	Assuming that Theorem \ref{thm: main1} fails, by Ramsey property and compactness (model theoretic jargon: extracting a generalized indiscernible) we can find some large probability space $(X,\mathcal{B}, \mu)$, $0<r<s<1$  and sets $A_{v,w} \in \mathcal{B}$ for all $v,w$ vertices of $\mathcal{OH}_2 =(E;V,W)$ so that:
	\begin{itemize}
	\item $(v,w) \in E \implies \mu(A_{v,w}) \geq s$ and $(v,w) \notin E \implies \mu(A_{v,w}) \leq r$;
		\item for any two isomorphic (as ordered bipartite graphs) substructures $H_1, H_2$ of $\mathcal{OH}_2$, 
		$$(\mathbf{1}_{A_{v,w}} : v,w \in H_1) =^{\textrm{dist}} (\mathbf{1}_{A_{v,w}} : v,w \in H_2).$$
	\end{itemize}
		
\subsection{(Relatively) Exchangeable random structures}\label{sec: exch arrays}

This indiscernibility guarantees certain ``exchangeability'' in the probabilistic sense. Exchangeable sequences (de Finetti, Section \ref{sec: de finetti}) and arrays (Aldous-Hoover-Kallenberg, see \cite{kallenberg2005probabilistic}) of random variables can be presented ``up to mixing'' using i.i.d.~random variables, and we need a certain generalization to relational structures which were studied recently by a number of authors \cite{crane2018relatively, ackerman2015representations, jahel2022invariant}.

\begin{definition}
(1) Let $\mathcal{L}' = \{R'_1, \ldots, R'_{k'}\}$, $R'_i$ a relation symbol of arity $r'_i$. By a \emph{random $\mathcal{L}'$-structure} we mean a (countable) collection of random variables 
	$ \left( \xi^{i}_{\bar{n}} : i \in [k'], \bar{n} \in \mathbb{N}^{r'_i} \right)$
	 on some probability space $\left(\Omega, \mathcal{F}, \mu \right)$ with $\xi^i_{\bar{n}}: \Omega \to \{0,1\}$. 
	 
(2) Let now $\mathcal{L} = \{R_1, \ldots, R_k\}$ be another relational language, with $R_i$ a relation symbol of arity $r_i$, and let $\mathcal{M} = (\mathbb{N}, \ldots)$ be a countable $\mathcal{L}$-structure with domain $\mathbb{N}$. We say that a random $\mathcal{L}'$-structure $ \left( \xi^{i}_{\bar{n}} : i \in [k'], \bar{n} \in \mathbb{N}^{r'_i} \right)$ is \emph{$\mathcal{M}$-exchangeable} if for any two finite subsets $A =\{a_1, \ldots, a_\ell\}, A' = \{a'_1, \ldots, a'_\ell\} \subseteq \mathbb{N}$ 
	\begin{gather*}
	\qftp_{\mathcal{L}}\left(a_1, \ldots, a_{\ell} \right) = \qftp_{\mathcal{L}}\left(a'_1, \ldots, a'_{\ell} \right) \implies \\
		\left( \xi^i_{\bar{n}} : i \in [k'], \bar{n} \in A^{r'_i} \right) =^{\textrm{dist}} \left( \xi^i_{\bar{n}} : i \in [k'], \bar{n} \in (A')^{r'_i} \right).
	\end{gather*}
\end{definition}

\subsection{A higher amalgamation condition on the indexing structure}
Let $\mathcal{K}$ be a collection of finite structures in a relational language $\mathcal{L}$ closed under isomorphism. 
\begin{definition}
	For $n \in \mathbb{N}_{\geq 1}$, we say that $\mathcal{K}$ satisfies the \emph{$n$-disjoint amalgamation property} (\emph{$n$-DAP}) if for every collection of $\mathcal{L}$-structures $\left( \mathcal{M}_i = (M_i, \ldots)  : i \in [n] \right)$ with $\mathcal{M}_i \in \mathcal{K}$, $M_i = [n] \setminus \{ i \}$ and
 $\mathcal{M}_i|_{[n] \setminus \{i,j\}} = \mathcal{M}_j|_{[n] \setminus \{i,j\}}$ for all $i \neq j \in [n]$,
	there exists an $\mathcal{L}$-structure $\mathcal{M} = (M, \ldots) \in \mathcal{K}$ such that $M = [n]$ and $\mathcal{M}|_{[n] \setminus \{ i \}} = \mathcal{M}_i$ for every $1 \leq i \leq n$.
\end{definition}
We say that an $\mathcal{L}$-structure $\mathcal{M}$ satisfies $n$-DAP if the collection of its finite substructures does. E.g., the generic $k$-hypergraph $\mathcal{H}_k$ satisfies $n$-DAP for all $n$ \cite[Proposition 9.6]{CheTow}, but $(\mathbb{Q}, <)$ fails $3$-DAP.

\subsection{Presentation for random relational structures}

\begin{fact}[Crane, Towsner \cite{crane2018relatively}; generalizing Aldous-Hoover-Kallenberg \cite{aldous1981representations, hoover1979relations, kallenberg1988spreading}]\label{fac: CraneTowsner} Let $\mathcal{L}' = \{R'_i : i \in [k']\}, \mathcal{L} = \{R_i : i \in [k]\}$ be finite relational languages with all $R'_i$ of arity at most $r'$, and $\mathcal{M} = (\mathbb{N}, \ldots)$ a countable homogeneous $\mathcal{L}$-structure that has $n$-DAP for all $n \geq 1$. Suppose that $ \left( \xi^{i}_{\bar{n}} : i \in [k'], \bar{n} \in \mathbb{N}^{r'_i} \right)$ is a random $\mathcal{L}'$-structure that is $\mathcal{M}$-exchangeable, such that the relations $R'_i$ are symmetric with probability $1$.

Then there exists a probability space  $(\Omega',\mathcal{F}',\mu')$, $\{0,1\}$-valued Borel functions $f_1, \ldots, f_{r'}$ and a collection of $\textrm{Uniform}[0,1]$ \emph{i.i.d.}~random variables $\left( \zeta_s : s \subseteq \mathbb{N}, |s| \leq r' \right)$ on $\Omega'$ so that
\begin{gather*}
	\left( \xi^i_{\bar{n}} : i \in [k'], \bar{n} \in \mathbb{N}^{r'_i} \right)=^{\textrm{dist}} 
	\left( f_i \left( \mathcal{M}|_{\range \bar{n}}, \left( \zeta_s \right)_{s \subseteq \range \bar{n}}\right) : i \in [k'], \bar{n} \in \mathbb{N}^{r'_i} \right),
\end{gather*}
where $\range \bar{n}$ is the set of its distinct elements, and $\subseteq$ denotes ``subsequence''.
\end{fact}

\subsection{Getting rid of the ordering}

Our counterexample from Section \ref{sec: exch counterex} is only guaranteed to be $\mathcal{OH}_{n}$-exchangeable (and the ordering is unavoidable in the Ramsey theorem for hypergraphs) --- but the presentation theorem in Fact \ref{fac: CraneTowsner} requires $n$-DAP (and linear orders fail $3$-DAP).
 However, using Fact \ref{fac: stability} inductively, we can show that $\mathcal{OH}_{n}$-exchangeability already implies $\mathcal{H}_{n}$-exchangeability (i.e., with respect to the reduct forgetting the ordering), using that the theory of probability algebras is \emph{stable}! (See \cite[Lemma 10.15]{CheTow} for the details.)
 
 Applying the exchangeable presentation to the counterexample, we finally reduce (modulo some mixing) to working with \emph{independent} random variables in the proof of Theorem \ref{thm: main1}, and can conclude the proof.

\section{Open questions and future directions}
\textbf{Question 1.} Our proof of Theorem \ref{thm: main1} is non-constructive and relies on a compactness argument. It would be interesting to obtain explicit bounds on $|G|$ and $\delta$ in terms of $H$ and $\varepsilon$. Do there exist infinitary/density versions of this result (similarly to Fact \ref{fac: basic})?

\noindent \textbf{Question 2.} Apart from \emph{$k$-partite} $k$-hypergraphs, which other classes of structures satisfy an analog of Theorem \ref{thm: main1}?
 E.g., there is a growing list of Ramsey classes of finite structures, for which also an appropriate analog of Fact \ref{fac: CraneTowsner} holds. The following example illustrates that these two properties alone are not sufficient:

\noindent \textbf{Example} (Tim Austin) Theorem \ref{thm: main1} does not hold for graphs (as opposed to \emph{bipartite} graphs). 
Indeed,  let $H$ be the triangle $K_3$, and $\varepsilon = 1/2$.  Consider any graph $G = (V,E)$.  On some probability space $(\Omega,\Sigma,\mu)$, let $(\pi_v: v \in V)$ be a process of independent uniform $\{0,1\}$-valued random variables, and consider the events $A_{v,w}$ defined by: 
$A_{v,w} := \left( {\pi_v \ne \pi_w} \right)$   if $(v,w) \in E$, and 
$A_{v,w} := \emptyset$   if $(v,w) \not\in E$.
Then $\mu(A_{v,w})$ is equal to $1/2$ if $(v,w) \in E$, but equal to $0$ if $(v,w) \not\in E$.  However, for any induced triangle in $G$, say with vertices $u,v,w$, we have
$\mu(A_{u,v} \cap A_{v,w} \cap A_{w,u}) = \mu(\pi_u \ne \pi_v \ne \pi_w \ne \pi_u) = \mu(\emptyset) = 0$.

As mentioned above, Theorem \ref{thm: main1} is the main ingredient in our proof that Keisler randomization of first-order theories preserves $n$-dependence,  for all $n \in \mathbb{N}_{\geq 1}$ (\cite[Corollary 11.3]{CheTow}).

\noindent \textbf{Question 3.} Apart from $n$-dependence, what are the other higher arity tameness notions from model theory preserved under Keisler randomization? E.g., is FOP$_n$ preserved?


\label{sect:bib}
\bibliographystyle{plain}
\bibliography{easychair}



\end{document}